\documentclass[11pt]{amsart}
\usepackage{amsmath}
\usepackage{amsthm}
\usepackage{amsfonts}
\usepackage{amssymb}
\usepackage[colorlinks=true]{hyperref}

\newtheorem{thm}{Theorem}[section]

\theoremstyle{definition}

\newcommand{\Q}{{\mathbb{Q}}}
\newcommand{\R}{{\mathbb{R}}}
\newcommand{\Z}{{\mathbb{Z}}}
\newcommand{\h}{{\mathrm{h}}}

\newcommand{\floor}[1]{\left\lfloor#1\right\rfloor}

\newcommand{\bp}{{b^\prime}}

\DeclareMathOperator*{\Rerm}{Re}
\DeclareMathOperator*{\Imrm}{Im}
\DeclareMathOperator*{\Log}{Log}

\let\lm=\lambda
\let\Lm=\Lambda
\let\ep=\epsilon
\let\abs=\envert
\let\sq=\sqrt

\newcommand{\acr}{\newline\indent}

\theoremstyle{remark}

\numberwithin{equation}{section}

\begin{document}
\title[The argument of an algebraic power]{A note on Laurent's paper on linear forms in two logarithms: the argument of an algebraic power}
\author[Tomohiro Yamada]{Tomohiro Yamada*}
\address{\llap{*\,}Center for Japanese language and culture,
                   The University of Osaka\acr
                   562-8678,
                   3-5-10, Semba-Higashi, Minoo, Osaka,
                   JAPAN}
\email{tyamada1093@gmail.com}


\begin{abstract}
In this note, we use Laurent's lower bound for linear forms in two logarithms in \cite{Lau2}
to give an improved lower bound for the argument of a power of a given algebraic number
which has absolute value one but is not a root of unity.
\end{abstract}

\maketitle

\section{Introduction}\label{intro}

Since Baker \cite{Bak1, Bak2} found lower bounds for linear forms in logarithms
\begin{equation}
b_1\log \alpha_1+b_2\log \alpha_2+\cdots +b_n\log \alpha_n
\end{equation}
with $\alpha_i$ complex algebraic numbers and $b_i$ integers,
many authors such as Matveev \cite{Mat} have given improved lower bounds for linear forms in logarithms of algebraic numbers.

Lower bounds for linear forms in two logarithms
\begin{equation}
\Lm=b_2\log \alpha_2-b_1\log \alpha_1,
\end{equation}
with $\alpha_1, \alpha_2$ two complex algebraic numbers and $b_1, b_2$ two positive integers
had already been given by Gel'fond \cite{Gel}
and several authors such as Laurent \cite{Lau1, Lau2} and Laurent, Mignotte and Nesterenko \cite{LMN}
have given improved lower bounds.

For any algebraic number $\alpha$ of degree $d$ over $\Q$,
we define the absolute logarithmic height of $\alpha$ by
\begin{equation}
\h(\alpha)=\frac{1}{d}\left(\log\abs{a}+\sum_{i=1}^d\log\max\{1, \abs{\alpha^{(i)}}\}\right),
\end{equation} 
where $a$ is the leading coefficient of the minimal polynomial of $\alpha$ over $\Z$
and $\alpha^{(i)} (i=1, \ldots, d)$ denote the conjugates of $\alpha$ in complex numbers.

As an application of their lower bound for linear forms in two logarithms,
Laurent, Mignotte and Nesterenko \cite[Th\'{e}or\`{e}me 3]{LMN} gave an lower bound
for this special logarithmic form
\begin{equation}
\Lm_0=b_2\log \alpha-b_1\pi i,
\end{equation}
where $\alpha$ is an algebraic number of absolute value one but not a root of unity,
$\log\alpha$ takes the principal branch, and  $b_1, b_2$ are positive integers.
Putting
\begin{equation}
\begin{split}
D= ~ & [\Q(\alpha): \Q]/2, \\
a= ~ & \max \{20, 10.98\abs{\log\alpha}+2D\h(\alpha)\}, \\
h= ~ & \max\{17, \sqrt{D}/10, D(\log (b_1/2a+b_2/68.9)+2.35)+5.03\},
\end{split}
\end{equation}
we have
\begin{equation}
\abs{\Lm_0}\geq -8.87ah^2.
\end{equation}
We note that the quantity $h$ here is denoted by $H$ in \cite{LMN}.
We use $h$ in order to be consistent with the notation \cite{Lau2}.

Later, Laurent \cite{Lau2} obtained the stronger lower bound for general linear forms in two logarithms
in the following form:
\begin{thm}\label{thm1}
Let $\Lm=b_2\log \alpha_2-b_1\log \alpha_1$ be a linear form of two logarithms
with $b_1, b_2$ positive integers and $\alpha_1, \alpha_2$ complex algebraic numbers.
Put $D=[\Q(\alpha_1, \alpha_2): \Q]/[\R(\alpha_1, \alpha_2): \R]$.

Let $K$ be an integer $\geq 2$ and $L, R_1, R_2, S_1, S_2$ be positive integers.
Let $\rho$ and $\mu$ be real numbers with $\rho>1$ and $1/3\leq\mu\leq 1$.
Put
\begin{equation}
\begin{split}
R= ~ & R_1+R_2-1, S=S_1+S_2-1, N=KL, g=\frac{1}{4}-\frac{N}{12RS}, \\
\sigma= ~ & \frac{1+2\mu-\mu^2}{2}, b=\frac{(R-1)b_2+(S-1)b_1}{2}\left(\prod_{k=1}^{K-1} k!\right)^{-2/(K^2-K)}.
\end{split}
\end{equation}
Let $a_1, a_2$ be positive real numbers such that
\begin{equation}\label{eq210}
a_i\geq \rho\abs{\log \alpha_i}-\log\abs{\alpha_i}+2D\h(\alpha_i)
\end{equation}
for $i=1, 2$.
Assume that
\begin{equation}\label{eq220}
\begin{split}
& \#\{\alpha_1^r \alpha_2^s: 0\leq r<R_1, 0\leq s<S_1\}\geq L, \\
& \#\{rb_2+sb_1: 0\leq r<R_2, 0\leq s<S_2\}>(K-1)L
\end{split}
\end{equation}
and
\begin{equation}\label{eq221}
K(\sigma L-1)\log\rho-(D+1)\log N-D(K-1)\log b-gL(Ra_1+Sa_2)>\ep(N),
\end{equation}
where
$\ep(N)=2\log (N! N^{-N+1}(e^N+(e-1)^N))/N$.

Then $\abs{\Lm^\prime}>\rho^{-\mu KL}$,
where
\begin{equation}
\Lm^\prime=\Lm \max\left\{\frac{LS e^{LS\abs{\Lm}/(2b_2)}}{2b_2} , \frac{LR e^{LR\abs{\Lm}/(2b_1)}}{2b_1}\right\}.
\end{equation}
\end{thm}

However, Laurent has not given an improved lower bound for the special logarithmic form $\Lm_0$.
Among two-logarithmic forms, this special form may be of some interest and
improving lower bounds for this logarithmic form may have some applications.
The purpose of this note is to deduce an improved lower bound for the special logarithmic form $\Lm_0$
from Theorem \ref{thm1}.

\begin{thm}\label{thm2}
Let
\begin{equation}
\Lm_1=b_2\log \alpha-\frac{b_1\pi i}{2},
\end{equation}
where $b_1, b_2$ are nonzero integers and $\alpha$ is an complex algebraic number of absolute value one
but not a root of unity.
Put
\begin{equation}
\begin{split}
a= ~ & 9\pi+2D\h(\alpha), \quad
\bp=\frac{\abs{b_1}}{a}+\frac{\abs{b_2}}{9\pi}, \quad
D=\frac{[\Q(\alpha): \Q]}{2}, \\
h= ~ & \max\{17, D, D(\log \bp+2.96)+0.01\}.
\end{split}
\end{equation}
Then,
\begin{equation}\label{eq231}
\log\abs{\Lm_1}>-(2.76701a+0.12945)h^2.
\end{equation}
Moreover, taking
$h^\prime=\max\{1000, D, D(\log \bp+2.96)+0.01\}$, we have
\begin{equation}\label{eq232}
\log\abs{\Lm_1}>-(1.86151a+0.00143)h^{\prime \, 2}.
\end{equation}
\end{thm}

Recently, combining Theorem \ref{thm2} and three-logarithmic forms,
the author \cite{Ymd} proved that there exist only finitely many three-term Machin-type formulae
for integer multiples of $\pi/4$ which are non-degenerate (i.e. not derived from two-term Machin-type formulae)
and gave explicit upper bounds for the sizes of the variables.

It immediately follows from Theorem \ref{thm2} that
if $\alpha$ is a complex algebraic number of absolute value one and positive real part but not a root of unity,
$n$ is a nonzero integer, then,
using the same notation as in Theorem \ref{thm2} with
$b_1$ the nearest integer to $2n\abs{\arg\alpha}/\pi$ and $b_2=n$,
\begin{equation}
\log \abs{\arg (\alpha^n)}>-\min \{(2.76701a+0.12945)h^2, (1.86151a+0.00143)h^{\prime \, 2}\},
\end{equation}
where the argument $\arg z$ is taken from the range $-\pi<\arg z\leq \pi$.

We used PARI-GP \cite{PARI} for our calculations.
Our script can be downloaded from \\
\url{https://drive.google.com/file/d/1HYTuCuveIHwgs84hiy8R4_qBacnkBBP0/}.

\section{Preliminaries to the proof}
We note that we work with a slightly generalized form $\Lm_1$ rather than $\Lm_0$.

We begin by noting that we may assume that $b_2>0$.
Moreover, we may assume that $\alpha$ has positive imaginary part
by changing $\alpha$ into $\overline\alpha$ if necessary.

For now, we limit ourselves that $\log\alpha$ takes its principal value.
Taking the assumption made above, we have $0<\Imrm\log\alpha<\pi$.
We see that $b_1, b_2>0$ since if $b_1$ and $b_2$ have opposite signs, then $\abs{\Lm_1}\geq \pi/2$ and
Theorem \ref{thm2} immediately follows.
Moreover, we limit ourselves in the case $\Rerm \alpha>0$ and therefore
we have $0<\Imrm\log\alpha<\pi/2$ and $b_1, b_2$ are positive integers.
Indeed, it is easy to prove Theorem \ref{thm2} in the general case
once we prove Theorem \ref{thm2} for $\Rerm\alpha>0$ and $0<\Imrm\log\alpha<\pi/2$
as we shall see later.

If $d=\gcd(b_1, b_2)>1$, then we divide $b_i$'s by $d$
to have another logarithmic form $\Lm_1/d=(b_2/d)\log\alpha-(b_1/d) \pi i/2$.
If Theorem \ref{thm2} holds for $\abs{\Lm_1/d}$, then this would give
the desired lower bound for $\abs{\Lm_1}$.
Thus we may assume that $\gcd(b_1, b_2)=1$.

Moreover, we may assume that $\bp>5.52h^2$.
Indeed, if $\bp\leq 5.52h^2$,
then, after observing that $ah>153.85\pi>480$ and $D\log 2<h\log 2<0.002ah^2$, Liouville's inequality
in the form of Exercise 3.6.b in p. 109 of \cite{Wal} gives
\begin{equation}
\log \abs{\Lm_1}\geq -\bp D\h(\alpha) -D\log 2>-2.76ah^2-D\log 2>-2.762ah^2.
\end{equation}

We set 
\begin{equation}\label{eq30}
\begin{split}
\mu= ~ & 0.59, \quad \rho=18, \quad \alpha_1=i, \quad \alpha_2=\alpha, \\
a_1= ~ & \frac{\rho\pi}{2}, \quad a_2=\frac{\rho\pi}{2}+2D\h(\alpha)=a\geq a_1.
\end{split}
\end{equation}
We see that $\alpha_i$'s and $a_i$'s satisfy the condition \eqref{eq210} in Theorem \ref{thm1}
since we have assumed that $\log\alpha=\kappa i$ with $0<\kappa<\pi/2$.

According to p.336, (8) and succeeding formulae of \cite{Lau2}, we put
\begin{equation}\label{eq31}
\begin{split}
\sigma & =\frac{1+2\mu-\mu^2}{2}=0.91595, \quad \lm=\sigma\log\rho=2.647436\cdots, \\
H & =\frac{h}{\lm}+\frac{1}{\sigma}, \quad L_0=H+\sqrt{H^2+\frac{1}{4}}, \quad L=\floor{L_0+\frac{1}{2}},
\end{split}
\end{equation}
and $k=(v_1(L)+\sqrt{v_1(L)^2+4v_0(L)v_2(L)})^2/(2v_2(L))^2$ to be the positive real number such that
$\sq{k}$ satisfies the quadratic equation
$v_2(L) k-v_1(L)\sqrt{k}-v_0(L)=0$,
where
\begin{equation}
v_0(x)=\frac{1}{4a_1}+\frac{4}{3a_2}+\frac{x}{12a_1}, \quad
v_1(x)=\frac{x}{3}, \quad
v_2(x)=\lm (x-H).
\end{equation}
We note that $v_0(x), v_1(x), v_2(x)$ correspond to $W, V, U$ in p. 336, l.5 of \cite{Lau2} respectively.
Moreover, as in (7.1) of \cite{LMN} and p.336, l.10-11 of \cite{Lau2}, we set
\begin{equation}\label{eq32}
\begin{split}
K= ~ & 1+\floor{kLa_1 a_2}, \quad R_1=4, \quad S_1=\floor{\frac{L+3}{4}}, \\
R_2= ~ & 1+\floor{\sq{(K-1)La_2/a_1}}, \quad S_2=1+\floor{\sq{(K-1)La_1/a_2}}.
\end{split}
\end{equation}

We see that $L^\pm=L_0\pm 1/2$ satisfies the quadratic equation $(L^\pm)^2-2L_0(L^\pm-H)=0$
and therefore $(L^\pm)^2/(L^\pm-H)=2L_0$.
Moreover, we can easily see that $L^-<2H<L^+$ and the function $x^2/(x-H)$ takes its minimum at $x=2H$ and is monotone below and above $x=2H$.
Hence, we obtain
\begin{equation}\label{eq331}
\frac{L^2}{L-H}\leq \frac{(L^{\pm})^2}{L^\pm-H}=2L_0
\end{equation}
and
\begin{equation}\label{eq332}
\sq{k}>\frac{v_1(L)}{v_2(L)}=\frac{L}{3\lm(L-H)}>\frac{L^+}{3\lm(L^+ -H)}>0.2437.
\end{equation}
Since both $v_0(x)/v_2(x)$ and $v_1(x)/v_2(x)$ are clearly decreasing for $x>H$, we have
\begin{equation}\label{eq332b}
\sq{k}<\frac{v_1(L^-)}{2v_2(L^-)}+\sqrt{\left(\frac{v_1(L^-)}{2v_2(L^-)}\right)^2+\frac{v_0(L^-)}{v_2(L^-)}}
<0.2795.
\end{equation}
Moreover, we have
\begin{equation}\label{eq333}
H>7.5, 15\leq L\leq L_0+\frac{1}{2}<0.92h, K>kLa_1 a_2>712.
\end{equation}
Indeed, since $h\geq 17$, \eqref{eq31} immediately gives $H>7.5$, $L_0>15$, and $L\geq 15$.
We observe that $(L_0+1/2)/h$ is monotonically decreasing for $h\geq 17$.
Hence, $(L_0+1/2)/h<0.92$ for $h\geq 17$.
Combining the fact that $L\geq 15$ with \eqref{eq30} and \eqref{eq332}, we have
$kLa_1 a_2>15(0.2437\times 9\pi)^2>712$.

We see that
\begin{equation}
\sqrt{k}L=\frac{L^2}{6U}+\frac{1}{2}\sqrt{\left(\frac{L^2}{3U}\right)^2+\frac{L^3}{3a_1 U}+\left(\frac{1}{a_1}+\frac{16}{3a_2}\right)\frac{L^2}{U}},
\end{equation}
where $U=\lm(L-H)$.
Using \eqref{eq331} and then recalling from \eqref{eq30} that $a_2\geq a_1$, we have
\begin{equation}\label{eq334}
\sqrt{k}L\leq \frac{L_0}{3\lm}
+\sqrt{\left(\frac{L_0}{3\lm}\right)^2+\frac{2L_0}{\lm}\left(\frac{19+L^+}{12a_1}\right)}<0.239537h.
\end{equation}

\section{Confirmation of the conditions of Theorem \ref{thm1}}

In this section, we shall confirm the conditions of Theorem \ref{thm1}.

In order to obtain an upper bound for $gL(Ra_1+Sa_2)$,
we follow the proof of Lemme 9 of \cite{LMN}.
We begin by quoting the upper bound
\begin{equation}\label{eq341}
gL(Ra_1+Sa_2)\leq \frac{L}{4}(R_1 a_1+S_1 a_2)+\frac{L^{3/2}\sq{(K-1)a_1 a_2}}{2}-\frac{KL^2}{12}\left(\frac{a_1}{S}+\frac{a_2}{R}\right)
\end{equation}
from (5.19) of \cite{LMN}.

As in \cite{LMN}, using the identity $\frac{1}{x+y}=\frac{1}{x}-\frac{y}{x^2}+\frac{y^2}{x^2(x+y)}$,
we obtain
\begin{equation}
\begin{split}
\frac{1}{R}>\frac{1}{R_1+\sq{(K-1)La_2/a_1}}
= & ~ \frac{1}{\sqrt{(K-1)La_2/a_1}}-\frac{R_1}{(K-1)La_2/a_1} \\
& +\frac{a_1 R_1^2}{(K-1)La_2(R_1+\sqrt{(K-1)La_2/a_1})}
\end{split}
\end{equation}
and
\begin{equation}
\begin{split}
\frac{1}{S}>\frac{1}{S_1+\sq{(K-1)La_1/a_2}}
= & ~ \frac{1}{\sqrt{(K-1)La_1/a_2}}-\frac{S_1}{(K-1)La_1/a_2} \\
& +\frac{a_2 S_1^2}{(K-1)La_1(S_1+\sqrt{(K-1)La_1/a_2})}.
\end{split}
\end{equation}
These lower bounds yield that
\begin{equation}
\begin{split}
& KL^2\left(\frac{a_1}{S}+\frac{a_2}{R}\right)>(K-1)L^2\left(\frac{a_1}{S}+\frac{a_2}{R}\right) \\
& >2L^{3/2}\sqrt{(K-1)a_1a_2}-L(R_1 a_1+S_1 a_2)+\frac{a_2 L S_1^2}{S_1+\sqrt{(K-1)La_1/a_2}} \\
& \qquad +\frac{a_1 LR_1^2}{R_1+\sqrt{(K-1)La_2/a_1}}.
\end{split}
\end{equation}
Now, \eqref{eq341} gives
\begin{equation}\label{eq342}
\begin{split}
gL(Ra_1+Sa_2)< & ~ \frac{L}{3}(R_1 a_1+S_1 a_2)+\frac{L^{3/2}\sq{(K-1)a_1 a_2}}{3} \\
& -\frac{a_2 L S_1^2}{12(S_1+\sqrt{(K-1)La_1/a_2})} \\
& -\frac{a_1 L R_1^2}{12(R_1+\sqrt{(K-1)La_2/a_1})}.
\end{split}
\end{equation}
Recalling that $R_1=4, S_1=\floor{(L+3)/4}\geq L/4$ and $K-1<kLa_1 a_2$, we have
\begin{equation}\label{eq343}
\begin{split}
gL(Ra_1+Sa_2)< & ~ \frac{L}{3}\left(4a_1+\frac{a_2(L+3)}{4}\right)+\frac{\sqrt {k} L^2 a_1 a_2}{3} \\
& -\frac{a_2 L^2}{48+192a_1\sqrt{k}}-\frac{4a_1 L}{12+3a_2 L\sqrt{k}} \\
< & ~ \left(\frac{\sqrt{k}}{3}+\frac{1}{12a_1}\right)a_1 a_2 L^2+\left(\frac{4}{3}a_1+\frac{a_2}{4}\right)L.
\end{split}
\end{equation}

Now we follow the proof of Lemme 10 of \cite{LMN}.
We put $\delta_1=0.044$.
From \eqref{eq332} and \eqref{eq333}, we see that
\begin{equation}
\frac{R_1-1}{R_2-1}<\frac{3}{\sq{(K-1)La_2/a_1}-1}<0.03<\delta_1
\end{equation}
and
\begin{equation}
\frac{S_1-1}{S_2-1}<\frac{S_1}{S_2}<\frac{1+3/L}{4a_1\sq{k}}\sq{\frac{K}{K-1}}< 0.044=\delta_1.
\end{equation}
Hence, we have
\begin{equation}
\log b<\log\bp+\frac{3}{2}+\log\left(\frac{1+\delta_1}{2\sq{k}}\right)+f_1(K)-\frac{\log(2\pi K/\sq{e})}{K-1},
\end{equation}
where
\begin{equation}
f_1(x)=\frac{1}{2}\log\left(\frac{x}{x-1}\right)+\frac{\log x}{6x(x-1)}+\frac{\log(x/(x-1))}{x-1}.
\end{equation}
\eqref{eq333} implies that $f_1(K)<f_1(712)<0.00072$.
Moreover, it follows from \eqref{eq332} that
$f_2(K):=f_1(K)+3/2+\log((1+\delta_1)/2\sq{k})<2.96$ and therefore
\begin{equation}\label{eq345}
\log b<\frac{h-\delta_2}{D}-\frac{\log(2\pi K/\sq{e})}{K-1},
\end{equation}
where we put $\delta_2=0.01$.

From \eqref{eq342} and \eqref{eq345}, we see that the left of \eqref{eq221} is at least
\begin{equation}
\begin{split}
& KL\lm-K\log\rho-(D+1)\log(KL)-(K-1)(h-\delta_2)+D\log(2\pi K/\sq{e}) \\
& ~ -\left(\frac{\sqrt{k}}{3}+\frac{1}{12a_1}\right)a_1 a_2 L^2+\left(\frac{4}{3}a_1+\frac{a_2}{4}\right)L \\
& > \left(L\left(k\lm-\frac{\sqrt{k}}{3}-\frac{1}{12a_1}\right)-k\lm H-\frac{1}{4a_1}-\frac{4}{3a_2}\right)L a_1 a_2 \\
& ~ +\delta_2(K-1)+h+D\log(2\pi K/\sq{e})-(D+1)\log (KL) \\
& = \Phi L a_1 a_2+\Theta,
\end{split}
\end{equation}
say.
We can easily see that $\Phi=v_2(L) k-v_1(L)\sqrt{k}-v_0(L)=0$ and therefore
\begin{equation}
K(\sigma L-1)\log\rho-(D+1)\log N-D(K-1)\log b-gL(Ra_1+Sa_2)>\Theta.
\end{equation}

Now we would like to show that $\Theta>\ep(N)$.
Our argument is similar to the argument in Section 3.2.2 of \cite{Lau2}.
Observing that $h-\delta_2>D(f_2(K)+\log b^\prime)$, we have $\Theta\geq (D-1)\Theta_0+\Theta_1$, where
\begin{equation}
\begin{split}
\Theta_0= ~ & \log\bp+f_2(K)-\log L+\log(2\pi/\sq{e}), \\
\Theta_1= ~ & \delta_2 K-\log K-2\log L+\log \bp+f_2(K)+\log(2\pi/\sq{e}).
\end{split}
\end{equation}
We recall the assumption that $\bp>5.52h^2$
and we see that $L\leq L^+<h$ from \eqref{eq333}.
Thus we obtain
\begin{equation}
\Theta_0>\log (5.52h)+f_2(K)+\log(2\pi/\sq{e})>0
\end{equation}
and
\begin{equation}
\Theta_1>\log 4+\delta_2 K-\log K+f_2(K)+\log(2\pi/\sq{e})>\delta_2 K-\log K>0.004.
\end{equation}
On the other hand, \eqref{eq333} gives that $N=KL>10000$ and, using Stirling's formula in the form
given in Section II.9 of \cite{Fel} or \cite{Rob}
we have $\ep(N)<\ep(10000)<0.003$.
This implies that our values for $k, L, R_1, S_1, R_2, S_2$ satisfy \eqref{eq221}.

Now we shall confirm \eqref{eq220}.
Since $\alpha_2=\alpha$ is not a root of unity,
$\alpha_1^r \alpha_2^s$ $(0\leq r\leq 3, 0\leq s\leq S_1-1)$ take $4S_1\geq L$ different values
and therefore the former part of \eqref{eq220} holds.

It follows from \eqref{eq334} that
$R_2-1<\sqrt{(K-1)La_2/a_1}<\sqrt{k}La_2<a_2 h/4<2a_2 h^2$
and, similarly, $S_2-1<a_1 h/4<2 a_1 h^2$.
Since we have assumed that $\bp>5.52h^2$,
$b_1>2a_2 h^2>R_2-1$ or $b_2>2a_1 h^2>S_2-1$.

Thus we can see that $R_2-1<b_1$ or $S_2-1<b_2$.
If we have $r_1 b_2-s_1 b_1=r_2 b_2-s_2 b_1$
for some integers $r_1, r_2, s_1, s_2$ with $0\leq r_1, r_2\leq R_2-1, 0\leq s_1, s_2\leq S_2-1$,
then $(r_1-r_2) b_2=(s_1-s_2)b_1$ and $\abs{r_1-r_2}\leq R_2-1, \abs{s_1-s_2}\leq S_2-1$.
Since we have assumed that $\gcd(b_1, b_2)=1$,
$r_1\equiv r_2\pmod{b_1}$ and $s_1\equiv s_2\pmod{b_2}$.
If $R_2-1<b_1$, then $r_1=r_2$.
If $S_2-1<b_2$, then $s_1=s_2$.
Hence, we must have $r_1=r_2$ and $s_1=s_2$.
This yields that $r b_2-s b_1$ $(0\leq r\leq R_2-1, 0\leq s\leq S_2-1)$ take $R_2 S_2>(K-1)L$ different values.
Hence, the latter part of \eqref{eq220} also holds.

Thus we have confirmed that Theorem \ref{thm1} holds with our choice of parameters.

\section{Computation of the constants}

Now we apply Theorem \ref{thm1} to obtain
$\log \abs{\Lm_1^\prime}>-\mu KL\log \rho$,
where
\begin{equation}
\Lm_1^\prime=\Lm_1 \max\left\{\frac{LS e^{LS\abs{\Lm_1}/(2b_2)}}{2b_2} , \frac{LR e^{LR\abs{\Lm_1}/(2b_1)}}{2b_1}\right\}.
\end{equation}

By \eqref{eq332} and \eqref{eq334}, we have
\begin{equation}\label{eq40}
KL<L(1+kLa_1 a_2)<kL^2\left(a_1 a_2+\frac{1}{kL}\right)<0.057378(a_1 a_2+1.12253)h^2.
\end{equation}
Thus, recalling that $a_1=\rho\pi/2$ and $a_2=a$ from \eqref{eq30},
we obtain $KL<(1.62233a+0.064409)h^2$.
Combined with the values of $\mu$ and $\rho$ in \eqref{eq30} again,
we obtain $\mu KL\log\rho<(2.7666a+1.0984)h^2$.
Now Theorem \ref{thm1} gives
\begin{equation}\label{eq41}
\log\abs{\Lm_1^\prime}>-(2.7666a+1.0984)h^2>-2.7705ah^2.
\end{equation}

We may assume that $\log \abs{\Lm_1}<-(2.7666a+1.0984)h^2$.
We see that
\begin{equation}
R<\frac{L+3}{4}+\sqrt{k} La_2<0.75+0.23h+0.239537ha<0.25ah
\end{equation}
and
$S<4+\sqrt{k} L a_1<4+0.239537h a_1<0.248ah$.
Thus, we obtain $LR, LS<0.23ah^2$.

We observe that $\log(x)/x$ is monotonically decreasing for $x>e$ and
\begin{equation}\label{eq42}
\log (ah^2)=\log a+\log(h^2)<(0.00041a+0.01961)h^2<0.00111ah^2.
\end{equation}
Hence, we see that $\log \max\{LR, LS\}+\log \abs{\Lm_1}<-2.76939h^2$ and
\begin{equation}
\begin{split}
\max\{LR\abs{\Lm_1}+\log(LR), LS\abs{\Lm_1}+\log(LS)\}< & ~ e^{-2.76939h^2}+\log(0.23ah^2) \\
< & ~ (0.00041a+0.01961)h^2.
\end{split}
\end{equation}
This immediately gives that
\begin{equation}
\log\abs{\Lm_1}>\log\abs{\Lm_1^\prime}-(0.00041a+0.01961)ah^2>-(2.76701a+0.12945)h^2.
\end{equation}
This proves \eqref{eq231} provided that $\Rerm\alpha>0$.

We can prove \eqref{eq232} in a quite similar way.
We take $\rho=22.5$ and $\mu=0.62$.
From \eqref{eq31} with $h$ replaced by $h^\prime\geq 1000$,
we see that $H>347.251$, $L_0>694.5$, $L\geq \floor{L_0+1/2}\geq 695$, and $\sqrt{k}>0.2306$.
Hence, we obtain
\begin{equation}
K>46164, (R_1-1)/(R_2-1)<0.001, (S_1-1)/(S_2-1)<0.031
\end{equation}
and we can confirm the condition of Theorem \ref{thm1}.
We immediately have $1/(kL)<0.02706$ and,
using the fact that $\sqrt{k}L<0.1651804h^\prime$ instead of \eqref{eq334},
\eqref{eq40} becomes $KL<0.1651804^2(a_1 a_2+0.02706)h^2$ and
\eqref{eq41} becomes $\log\abs{\Lm_1^\prime}>-(1.8615a+0.00143)h^{\prime \, 2}$ now.
Observing that $\log(ah^{\prime \, 2})<0.00001ah^{\prime \, 2}$ in place of \eqref{eq42},
we have \eqref{eq232} for $\alpha$ with a positive real part.

We must prove Theorem \ref{thm2} also for $\alpha$ with a negative real part.
Like above, we may assume that $\alpha$ has a positive imaginary part and $b_1, b_2>0$.
Hence, we have $\pi/2<\arg\alpha<\pi$.
Then, we can apply Theorem \ref{thm2} to
\begin{equation}
-\Lm_1=-b_2\log \alpha+\frac{b_1\pi i}{2}=b_2\log(-\alpha)-\frac{(2b_2-b_1)\pi i}{2}
\end{equation}
and obtain \eqref{eq231} and \eqref{eq232} with $b_1$ replaced by $2b_2-b_1$.

Thus we see that if $b_2\leq b_1$, then $\abs{2b_2-b_1}\leq b_1$ and \eqref{eq231} and \eqref{eq232} hold.
If $b_1<b_2$, then we must have $b_2\geq b_1+1$.
However, since $\pi/2<\arg\alpha<\pi$, we have
\begin{equation}
\abs{\Lm_1}=\abs{b_2\arg \alpha-\frac{b_1\pi}{2}}>\frac{(b_2-b_1)\pi}{2}\geq\frac{\pi}{2}
\end{equation}
and \eqref{eq231} and \eqref{eq232} clearly hold.

Finally, we must settle the case $\log\alpha$ can take principal value.
Writing $\log\alpha=2\pi k i+\Log\alpha$ for some integer $k$, where $\Log\alpha$ denotes the principal
value of $\log\alpha$, we have
\begin{equation}
\Lm_1=b_2\log \alpha-\frac{b_1\pi i}{2}=b_2\Log\alpha-\frac{(b_1-4k)\pi i}{2}
=b_2\Log\alpha-\frac{b_1^\prime \pi i}{2}
\end{equation}
by taking $b_1^\prime=b_1-4k$.
We are reminded of the assumption that $\Imrm\alpha>0$ to see that $0<\Imrm\Log\alpha<\pi$.
Moreover, we may assume that $b_1^\prime\geq 0$ since otherwise $\abs{\Lm_1}=\Imrm\Lm_1\geq \pi/2$
recalling that $b_2>0$.

If $k\geq 0$, then $\Imrm\log\alpha>0$ and therefore $b_1$ and $b_2$ must have the same sign.
Hence, we have $0\leq b_1^\prime=b_1-4k\leq b_1$.
Now the theorem follows from the principal case.

If $k<0$, then $\Imrm\log\alpha\leq \Imrm\Log\alpha-2\pi<-\pi$.
Hence, $b_1$ and $b_2$ must have opposite signs and $b_1<0$, recalling that $b_2>0$ again.
On the other hand, we see that $b_1^\prime\geq 0$ 
If $b_1^\prime>\abs{b_1}$, then, observing that $\abs{\log\alpha}>\pi>\Log\alpha$, we have
\begin{equation}
b_2\Log\alpha<\abs{b_2\log \alpha}\leq \frac{\abs{b_1}\pi}{2}+\abs{\Lm_1}\leq\frac{(b_1^\prime-1)\pi}{2}+\abs{\Lm_1}
\end{equation}
and $\abs{\Lm_1}\geq \pi/2$.
If $0\leq b_1^\prime\leq \abs{b_1}$, then the theorem immediately follows from the principal case.
This completes the proof of Theorem \ref{thm2}.

\section{Acknowledgements}
We would like to thank the anonymous referee for very careful reading and many useful comments.

{}
\end{document}